# Special Functions Associated with Root Systems: Recent Progress

Tom H. Koornwinder

*University of Amsterdam, Department of Mathematics,*
*Plantage Muidergracht 24, 1018 TV Amsterdam, The Netherlands*
*e-mail:* `thk@fwi,uva.nl`

**1. Introduction**

In this paper I present a brief survey of the active area of Special Functions associated with Root Systems. The article is intended for a general mathematical audience. It will not suppose prerequisites on either special functions or root systems. It will also skip many technical details. Some early work in this area (root systems $BC_2$ and $A_2$) was done by the author [17] and Sprinkhuizen [32]. During the last ten years important break-throughs were made by Heckman and Opdam [11], [27], [28], [12], Macdonald [23], [24], Dunkl [8] and Cherednik [3], [4], [5].

A lot of the motivation for the subject of this paper comes from analysis on semisimple Lie groups. Spherical functions on Riemannian symmetric spaces of the compact or non-compact type can be written as special functions depending on parameters which assume only special discrete values. In the one-variable cases these special functions are classical, also for parameter values without group theoretic interpretation, but in the more-variable cases they were new. In the case of group theoretic interpretation, many properties of these special functions, as well as associated harmonic analysis, immediately follow by group theoretic arguments. The case of more general parameter values yields the special functions associated with root systems. Properties derived in the group case can still be formulated in the general case, but now as conjectures rather than theorems. This paper describes some of the progress which has been made in proving these conjectures. For convenience, I will restrict to the polynomial (compact) case, with Bessel functions as a sole exception. (For an introduction to the non-polynomial case see [29], [13].) Neither will I discuss the recent work on commuting operators with elliptic functions as coefficients. An important aspect of the whole theory, which will not be discussed very much in this paper, is the connection with completely integrable systems, for instance the generalized Calogero-Moser system.

Special functions associated with root systems have also been developed in the $q$-case, where $q$ is a deformation parameter giving back the earlier cases when $q = 1$. Motivation and development of the theory in the $q$-case has been quite different from the $q = 1$ case. Except for the case of Hall polynomials [22], theory was developed [23], [24], [18] without interpretation in group theory. But afterwards quantum groups looked very promising as

---

a natural setting for these polynomials. This had turned out to be true in the one-variable case [19], and very recently some interpretations of more-variable cases on quantum groups were found [25], [10].

In any case, a quantum group interpretation for generic values of the parameters cannot be expected. But, by Cherednik's work [3], [4], [5] we know already another algebraic setting for special functions associated with root systems: affine and graded Hecke algebras [20]. As shown by work of Opdam [29], this new algebraic context also allows harmonic analysis.

## 2. The one-variable case

In this section I will introduce three classical families of special functions, each depending on a real parameter $k \geq 0$, and such that the cases $k = 0$ and $k = 1$ are elementary. The three families are connected with each other by limit transitions. Later, for each of the families I will discuss generalizations which are associated with root systems.

**2.1. Bessel functions.** Consider *Bessel functions* in a non-standard notation:

$$\mathcal{J}_k(x) := \sum_{j=0}^{\infty} \frac{(-\frac{1}{4}x^2)^j}{(k+\frac{1}{2})_j \, j!} \quad (x \in \mathbb{R}). \tag{2.1}$$

Here we use the notation for *shifted factorial*:

$$(a)_j := a\,(a+1)\ldots(a+j-1) \quad (j = 1, 2, \ldots); \qquad (a)_0 := 1.$$

The function $\mathcal{J}_k$ is related to the Bessel function $J_\alpha$ in standard notation [9, Ch. VII] by

$$\mathcal{J}_k(x) = \frac{2^{k-\frac{1}{2}} \, \Gamma(k+\frac{1}{2})}{x^{k-\frac{1}{2}}} \, J_{k-\frac{1}{2}}(x).$$

Note that

$$\mathcal{J}_k(x) = \mathcal{J}_k(-x), \qquad \mathcal{J}_k(0) = 1. \tag{2.2}$$

The cases $k = 0$ and $k = 1$ yield elementary functions:

$$\mathcal{J}_0(x) = \cos x, \qquad \mathcal{J}_1(x) = \frac{\sin x}{x}. \tag{2.3}$$

The function $x \mapsto \mathcal{J}_k(\lambda x)$ ($\lambda \in \mathbb{R}$) is eigenfunction of a differential operator:

$$\left(\frac{d^2}{dx^2} + \frac{2k}{x}\frac{d}{dx}\right) \mathcal{J}_k(\lambda x) = -\lambda^2 \, \mathcal{J}_k(\lambda x).$$

It is the unique $C^\infty$ solution of this differential equation under conditions (2.2).

## 2.2. Ultraspherical polynomials.
Consider *ultraspherical* or *Gegenbauer polynomials* [9, §10.9], i.e. polynomials $C_n^k$ of degree $n$ on $\mathbb{R}$ such that

$$\int_0^\pi C_n^k(\cos x)\, C_m^k(\cos x)\, (\sin x)^{2k}\, dx = 0 \quad (n, m \in \mathbb{Z}_+,\ n \neq m).$$

Then the $C_n^k$ are determined up to a constant factor (in general, we will not use the standard normalization for Gegenbauer polynomials). For $k = 0, 1$ we have:

$$C_n^0(\cos x) = \text{const. } \cos(nx), \qquad C_n^1(\cos x) = \text{const. } \frac{\sin((n+1)x)}{\sin x}. \qquad (2.4)$$

The function $x \mapsto C_n^k(\cos x)$ is eigenfunction of a differential operator:

$$\left(\frac{d^2}{dx^2} + 2k \cot x \frac{d}{dx}\right) C_n^k(\cos x) = -n(n + 2k)\, C_n^k(\cos x).$$

For $(n_N)$ being a sequence of positive integers such that $n_N/N \to \lambda$ for some $\lambda \geq 0$ as $N \to \infty$, we have the limit result

$$\lim_{n \to \infty} \frac{C_{n_N}^k(\cos(x/N))}{C_{n_N}^k(1)} = \mathcal{J}_k(\lambda x).$$

## 2.3. $q$-Ultraspherical polynomials.
Let $0 < q < 1$ and define for any complex $a$:

$$(a; q)_\infty := \prod_{j=0}^\infty (1 - aq^j).$$

The infinite product converges because of the condition on $q$. We will consider *$q$-ultraspherical polynomials* [1] in a non-standard notation. These are polynomials $C_n^{k,q}$ of degree $n$ on $\mathbb{R}$ such that

$$\int_0^\pi C_n^{k,q}(\cos x)\, C_m^{k,q}(\cos x) \left|\frac{(e^{2ix}; q)_\infty}{(q^k e^{2ix}; q)_\infty}\right|^2 dx = 0 \quad (n, m \in \mathbb{Z}_+,\ n \neq m).$$

Then the $C_n^{k,q}$ are determined up to a constant factor. If we put

$$P_n(e^{ix}) := C_n^{k,q}(\cos x)$$

then $P_n$ is eigenfunction of a $q$-difference operator:

$$\frac{1 - q^k e^{2ix}}{1 - e^{2ix}} P_n(q^{\frac{1}{2}} e^{ix}) + \frac{1 - q^k e^{-2ix}}{1 - e^{-2ix}} P_n(q^{-\frac{1}{2}} e^{ix}) = (q^{-\frac{1}{2}n} + q^{\frac{1}{2}n+k})\, P_n(e^{ix}).$$

Note that the $P_n$ on the left hand side have arguments off the unit circle, while orthogonality is on the unit circle. The cases $k = 0$ and $k = 1$ are elementary as in (2.4) (not depending on $q$):

$$C_n^{0,q}(\cos x) = \text{const. } \cos(nx), \qquad C_n^{1,q}(\cos x) = \text{const. } \frac{\sin((n+1)x)}{\sin x}.$$

With suitable normalization there is the limit relation

$$\lim_{q \uparrow 1} C_n^{k,q}(\cos x) = C_n^k(\cos x).$$

The $q$-ultraspherical polynomials form a subclass of the *Askey-Wilson polynomials* [2]: a family of orthogonal polynomials depending, apart from $q$, on four non-trivial parameters.

### 2.4. Dunkl operators in one variable.
We will now generalize the elementary formulas

$$e^{i\lambda x} = \mathcal{J}_0(\lambda x) + i\lambda x\, \mathcal{J}_1(\lambda x) \quad \text{and} \quad \frac{d}{dx} e^{i\lambda x} = i\lambda\, e^{i\lambda x} \qquad (2.5)$$

(the first formula follows by (2.3)). Dunkl [8] generalized the operator $d/dx$ to a mixture of a differential and a reflection operator:

$$(D^{(k)} f)(x) := f'(x) + k\, \frac{f(x) - f(-x)}{x}. \qquad (2.6)$$

Note that this *Dunkl operator* sends smooth functions to smooth functions. Let us define a *generalized exponential function* in terms of Bessel functions (2.1) by

$$\mathcal{E}_k(\lambda x) := \mathcal{J}_k(\lambda x) + \frac{i\lambda x}{2k+1}\, \mathcal{J}_{k+1}(\lambda x). \qquad (2.7)$$

Then it follows immediately from well-known differential recurrence formulas for Bessel functions that

$$D^{(k)}\, \mathcal{E}_k(\lambda x) = i\lambda\, \mathcal{E}_k(\lambda x). \qquad (2.8)$$

Formulas (2.7) and (2.8) generalize the formulas in (2.5). The function $x \mapsto \mathcal{E}_k(\lambda x)$ is the unique $C^\infty$ function which equals 1 in 0 and which is eigenfunction with eigenvalue $i\lambda$ of $D^{(k)}$.

For $(D^{(k)})^2$ we compute

$$(D^{(k)})^2 f(x) = f''(x) + \frac{2k}{x} f'(x) - k\, \frac{f(x) - f(-x)}{x^2}.$$

Thus, on even functions $f$ the square of the Dunkl operator acts as the differential operator $(d/dx)^2 + 2k x^{-1}\, d/dx$. In particular, its action on

$$\mathcal{J}_k(\lambda x) = \tfrac{1}{2}(\mathcal{E}_k(\lambda x) + \mathcal{E}_k(-\lambda x))$$

yields

$$(D^{(k)})^2\, \mathcal{J}_k(\lambda x) = -\lambda^2\, \mathcal{J}_k(\lambda x).$$

## 3. Preliminaries about root systems

**3.1. Definition of root system.** Let $V$ be a $d$-dimensional real vector space with inner product $\langle .\,,.\rangle$. For $\alpha \in V\setminus\{0\}$ let $s_\alpha$ denote the orthogonal reflection with respect to the hyperplane orthogonal to $\alpha$ (cf. Fig. 1):

$$s_\alpha(\beta) := \beta - \frac{2\langle \beta, \alpha \rangle}{\langle \alpha, \alpha \rangle}\alpha \quad (\beta \in V).$$

A *root system* in $V$ (see [15]) is a finite subset $R$ of $V\setminus\{0\}$ which spans $V$ and which satisfies for all $\alpha, \beta \in R$ the two properties that

$$s_\alpha(\beta) \in R \quad \text{and} \quad \frac{2\langle \beta, \alpha \rangle}{\langle \alpha, \alpha \rangle} \in \mathbb{Z}.$$

Clearly, if $\alpha \in R$ then $-\alpha = s_\alpha(\alpha) \in R$. For convenience, we will restrict ourselves to the case of a *reduced root system*, i.e., a root system $R$ such that, if $\alpha, \beta \in R$ and $\alpha = c\beta$ for some $c \in \mathbb{R}$, then $c = \pm 1$. The so-called *irreducible* root systems can be classified as four infinite families $A_n, B_n, C_n, D_n$ of *classical root systems* and five *exceptional root systems* $G_2, F_4, E_6, E_7, E_8$. Here the subscript denotes the *rank* of the root system, i.e. the dimension of $V$. There is one infinite family of non-reduced irreducible root systems: of type $BC_n$.

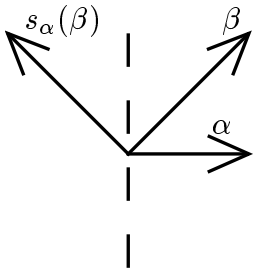 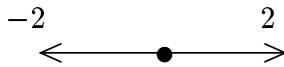 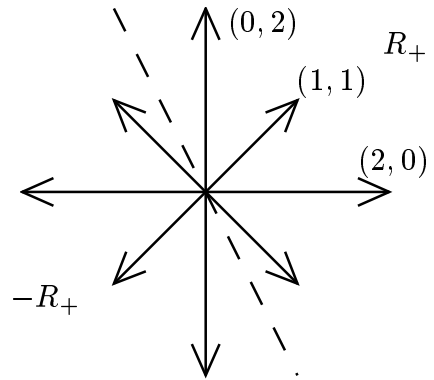

Fig. 1. Reflection $s_\alpha$  Fig. 2. Root system $A_1$  Fig. 3. Root system $C_2$

An example for $d = 1$ is the set $R := \{\pm 2\} \subset \mathbb{R}$ (root system of type $A_1$, cf. Fig. 2). An example for $d = 2$ is the set $R = R_+ \cup (-R_+)$, where $R_+ := \{(1,-1),(2,0),(1,1),(0,2)\} \subset \mathbb{R}^2$ (root system of type $C_2$, cf. Fig. 3). In general, when we have a root system $R$ in $V$ then we can write it as a disjoint union $R = R_+ \cup (-R_+)$, where $R_+$ and $-R_+$ are separated from each other by a hyperplane in $V$ through the origin. The choice for $R_+$ is not unique. The elements of $R$ are called *roots* and the elements of $R_+$ are called *positive roots*.

Let $GL(V)$ be the group of invertible linear transformations of $V$. The *Weyl group* $W$ of the root system $R$ is the subgroup of $GL(V)$ which is generated by the reflections $s_\alpha$ ($\alpha \in R$). The group $W$ is finite and it acts on $R$. It permutes the possible choices of $R_+$ in a simply transitive way.

## 3.2. Dunkl operators associated with $R$.

Let $R$ be a root system in $V$. Let $k: \alpha \mapsto k_\alpha: R \to [0, \infty)$ be a function which is $W$-invariant, i.e., which satisfies $k_{w\alpha} = k_\alpha$ for all $w \in W$ and all $\alpha \in R$. If $R$ is an irreducible (reduced) root system then the Weyl group is transitive on all roots of equal length and there are at most two different root lengths. Thus $k_\alpha$ then assumes at most two different values. See the above examples: one root length in $A_1$ and two root lengths in $C_2$. The function $k$ is called a *multiplicity function*. The reason for this name is that root systems naturally arise in the structure theory of real semisimple Lie algebras, where roots have an interpretation as joint eigenvalues of certain operators and the $k_\alpha$ then are (integer) multiplicities of such eigenvalues.

For $\xi \in V$ we will denote by $\partial_\xi$ the corresponding directional derivative. The *Dunkl operators* [8], [16] associated with the root system $R$ and the multiplicity function $k$ are defined as the operators $D_\xi^{(k)}: C^\infty(V) \to C^\infty(V)$ ($\xi \in V$) given by

$$(D_\xi^{(k)} f)(x) := (\partial_\xi f)(x) + \sum_{\alpha \in R_+} k_\alpha \langle \alpha, \xi \rangle \frac{f(x) - f(s_\alpha x)}{\langle \alpha, x \rangle}. \tag{3.1}$$

This definition is easily seen to be independent of the choice of $R_+$. In case of root system $A_1$ formula (3.1) reduces for $\xi := 1$ to formula (2.6). Note that the operator (3.1) consists of a term involving a first order derivative and terms involving reflection operators, just as we have seen in (2.6). It is an amazing fact, which can be proved in a rather straightforward way, that the operators $D_\xi^{(k)}$ commute:

$$[D_\xi^{(k)}, D_\eta^{(k)}] = 0 \quad (\xi, \eta \in V).$$

Let $\mathbb{D}^{(k)}$ be the algebra generated by the operators $D_\xi^{(k)}$. This is a commutative algebra. It can be shown that each $W$-invariant operator $D$ in $\mathbb{D}^{(k)}$, when restricted in its action to the $W$-invariant $C^\infty$ functions on $V$, coincides with a partial differential operator (so its reflection terms vanish when acting on a $W$-invariant function). The joint $W$-invariant eigenfunctions of the $W$-invariant operators in $\mathbb{D}^{(k)}$ are called *Bessel functions associated with $R$*. In the example $A_1$ things reduce to the one-variable considerations of §2.1 and §2.4. More generally, one may study the joint eigenfunctions of the full algebra $\mathbb{D}^{(k)}$ and one may try to do harmonic analysis for these eigenfunctions. A lot of satisfactory results have been obtained, see [16] and the references given there.

## 3.3. Weight lattice asscociated with $R$.

We still assume a root system $R$ in $V$. The *weight lattice* $P$ of $R$ is defined [15] by

$$P := \{ \lambda \in V \mid \frac{2\langle \lambda, \alpha \rangle}{\langle \alpha, \alpha \rangle} \in \mathbb{Z} \quad \text{for all } \alpha \in R \}.$$

The subset $P_+$ of *dominant weights* is then given by

$$P_+ := \{ \lambda \in P \mid \frac{2\langle \lambda, \alpha \rangle}{\langle \alpha, \alpha \rangle} \geq 0 \quad \text{for all } \alpha \in R_+ \}.$$

It is easily seen that $w(P) = P$ for $w \in W$, so the Weyl group acts on $P$. Moreover, it can be shown that each Weyl group orbit in $P$ has a one-point intersection with $P_+$:
$$\forall \lambda \in P \quad \mathrm{Card}(W\lambda \cap P_+) = 1.$$
Thus the dominant weights can be used as a set of representatives for the $W$-orbits in $P$.

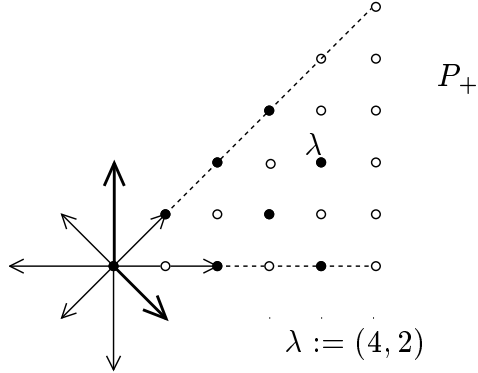

Fig. 4. Root system $C_2$ with dominant weights and the set $\{\mu \in P_+ \mid \mu \prec \lambda\}$

We introduce a partial ordering on $P$ which is induced by the root system: for $\lambda, \mu \in P$ we say that $\mu \prec \lambda$ if $\lambda - \mu = \sum_{\alpha \in R_+} m_\alpha \alpha$ for certain nonnegative integers $m_\alpha$. For root system $C_2$ the concepts of this subsection are illustrated in Fig. 4.

**3.4. Trigonometric polynomials associated with $R$.** Let $P$ be the weight lattice of a root system $R$ in $V$. For $\lambda \in P$ define the function $e^\lambda$ on $V$ by
$$e^\lambda(x) := e^{i\langle \lambda, x \rangle} \quad (x \in V).$$
Note that $e^\lambda e^\mu = e^{\lambda+\mu}$. Thus the space
$$\mathcal{A} := \mathrm{Span}\{e^\lambda \mid \lambda \in P\}$$
is an algebra: the *algebra of trigonometric functions on $V$* (with respect to $R$). For a function $f$ on $V$ write $(wf)(x) := f(w^{-1}x)$ ($w \in W$, $x \in V$). Then $we^\lambda = e^{w\lambda}$ ($w \in W$, $\lambda \in P$). Put
$$m_\lambda := \sum_{\mu \in W\lambda} e^\mu \quad (\lambda \in P_+).$$
Then the functions $m_\lambda$ are $W$-invariant and they form a basis of the space $\mathcal{A}^W$ of $W$-invariant elements in $\mathcal{A}$.

Let the *dual root lattice* $Q^{\vee}$ be defined by
$$Q^{\vee} := \{\lambda \in V \mid \langle \lambda, \mu \rangle \in \mathbb{Z} \quad \text{for all } \mu \in P\}.$$
This lattice gives rise to a torus
$$T := V/(2\pi Q^{\vee}).$$
Let $x \mapsto \dot{x}$ be the natural mapping of $V$ onto $T$. Then each function $f$ in $\mathcal{A}$ actually lives on $T$: $f(x) = \tilde{f}(\dot{x})$ for a suitable function $\tilde{f}$ on $T$.

In the example $A_1$ we have $P = \mathbb{Z}$, $P_+ = \{0, 1, 2, \ldots\}$, the algebra $\mathcal{A}$ is spanned by the functions $x \mapsto e^{inx}$ ($n \in \mathbb{Z}$) and the subalgebra $\mathcal{A}^W$ by the functions $1$ and $x \mapsto 2\cos(nx)$ ($n = 1, 2, \ldots$). The torus $T$ equals $\mathbb{R}/(2\pi\mathbb{Z})$.

## 4. Jacobi polynomials associated with $R$

**4.1. Definition of Jacobi polynomials for $R$.** Let $R$ be a root system in $V$ and let $k\colon R \to [0,\infty)$ be a $W$-invariant multiplicity function as before. Define a weight function $\delta_k$ on $T$ by

$$\delta_k(x) := \prod_{\alpha \in R_+} |2\sin(\langle \alpha, x \rangle)|^{2k_\alpha}. \tag{4.1}$$

This definition is independent of the choice of $R_+$. Define an inner product on the linear space $\mathcal{A}$ by

$$\langle f, g \rangle_k := \int_T f(x) \overline{g(x)} \delta_k(x) \, d\dot{x} \quad (f, g \in \mathcal{A}). \tag{4.2}$$

Here $d\dot{x}$ denotes Lebesgue measure on $T$, normalized such that the volume of $T$ is equal to 1.

The *Jacobi polynomial* $P_\lambda^{(k)}$ (cf. [11]) of "degree" $\lambda \in P_+$ and of "order" $k$ is an element of $\mathcal{A}^W$ of the form

$$P_\lambda^{(k)} = \sum_{\substack{\mu \in P_+ \\ \mu \prec \lambda}} c_{\lambda,\mu} \, m_\mu$$

such that $c_{\lambda,\lambda} = 1$ and

$$\langle P_\lambda^{(k)}, m_\mu \rangle_k = 0 \quad \text{if } \mu \in P_+ \text{ and } \mu \not\succeq \lambda. \tag{4.3}$$

Instead of (4.3) we can equivalently require that $P_\lambda^{(k)}$ satisfies the second order differential equation

$$\left( \Delta + \sum_{\alpha \in R_+} k_\alpha \cot(\tfrac{1}{2}\langle \alpha, x\rangle) \, \partial_\alpha \right) P_\lambda^{(k)}(x) = -\langle \lambda, \lambda + \textstyle\sum_{\alpha \in R_+} k_\alpha \alpha \rangle \, P_\lambda^{(k)}(x). \tag{4.4}$$

In the example $A_1$ we obtain that $P_n^{(k)}(x) = \text{const.}\, C_n^k(\cos x)$, where $C_n^k$ is the ultraspherical polynomial of §2.2. The case of the (non-reduced) root system $BC_1$ would have given us, more generally, the classical one-variable Jacobi polynomials.

**4.2. Three problems and their solutions.** As soon as the above definition of Jacobi polynomials associated with $R$ is given, three highly nontrivial questions can naturally be posed:

1. It follows immediately from the definition that the orthogonality

$$\langle P_\lambda^{(k)}, P_\mu^{(k)} \rangle_k = 0 \tag{4.5}$$

holds if $\mu \prec \lambda$ or $\lambda \prec \mu$. What about (4.5) if $\lambda$ and $\mu$ are not related in the partial ordering?

2. Prove the existence of a commutative algebra of differential operators with $d$ algebraically independent generators, such that the operators in this algebra have the $P_\lambda^{(k)}$ ($\lambda \in P_+$) as joint eigenfunctions. (Note that the operator in (4.4) can be taken as one of the generators.)

3. Give an explicit expression for $\langle P_\lambda^{(k)}, P_\lambda^{(k)} \rangle_k$, or rather for its two factors

$$\frac{\langle P_\lambda^{(k)}, P_\lambda^{(k)} \rangle_k}{\langle P_0^{(k)}, P_0^{(k)} \rangle_k} \quad \text{and} \quad \int_T \delta_k(x)\, d\dot{x}. \tag{4.6}$$

In the past few years all these questions have been answered in the positive sense. Let me give some indications.

- If problem 2 can be solved then the answer to 1 follows readily, cf. [11]. Indeed, we need sufficiently many differential operators having the $P_\lambda^{(k)}$ as eigenfunctions such that the joint eigenvalues, in their dependence on $\lambda$, separate the points of $P_+$.

- For certain special choices of $k$ the functions $P_\lambda^{(k)}$, renormalized such that $P_\lambda^{(k)}(0) = 1$, have an interpretation as spherical functions on compact symmetric spaces $G/K$, cf. [14]. (For instance, in case $A_1$ the ultraspherical polynomial $C_n^{\frac{1}{2}m-1}$ can be interpreted as spherical function on the $(m-1)$-dimensional sphere $SO(m)/SO(m-1)$.) Then problems 1, 2 and the first half of problem 3 can be solved by using the group theoretic interpretation. The orthogonality (4.5) for general $\lambda, \mu$ follows by Schur's orthogonality relations for matrix elements of irreducible unitary representations of $G$. The first expression in (4.6) was explicitly computed by Vretare [33] in terms of Harish-Chandra's $c$-function related to the spherical functions on the corresponding non-compact symmetric space. The algebra of differential operators in problem 2 can be obtained by taking the radial parts of the $G$-invariant differential operators on $G/K$.

- For the classical root systems question 2 could be answered in a positive way by giving explicit expressions for generators of the algebra, see [17] for $BC_2$ and $A_2$, and Olshanetsky & Perelomov [26], Sekiguchi [31] and Debiard [6] for the higher rank cases.

- Heckman and Opdam [11] have given positive answers to 2, and hence to 1, by use of deep transcendental arguments. This also solved part of Problem 3 (the first expression in (4.6)). In 1982 Macdonald [21] had already given conjectures for the explicit evaluation of the second expression in (4.6), which could be proved in a number of special cases.

- Problem 3 for general $\lambda$ was solved by Opdam [28] by using so-called *shift operators* [27]. The most simple example, for case $A_1$, of such operators is the following pair of differential recurrence relations for Gegenbauer polynomials:

$$\frac{d}{dx} C_n^k(x) = \text{const.}\, C_{n-1}^{k+1}(x),$$

$$\left( (1-x^2)^{-k+\frac{1}{2}} \frac{d}{dx} \circ (1-x^2)^{k+\frac{1}{2}} \right) C_{n-1}^{k+1}(x) = \text{const.}\, C_n^k(x).$$

By use of these two formulas we can write $\int_{-1}^{1}(C_n^k(x))^2\,(1-x^2)^{k-\frac{1}{2}}\,dx$ as an explicit constant times $\int_{-1}^{1}(C_{n-1}^{k+1}(x))^2\,(1-x^2)^{k+\frac{1}{2}}\,dx$. Opdam's shift operators in general have a similar structure of lowering $\lambda$ and raising $k$, or conversely. The case of root system $BC_2$ was already considered in [17], [32].

**4.3. Dunkl type operators.** Some years after Heckman first solved the problems 1 and 2 of the previous subsection he discovered a dramatical simplification [12] for proving these results. For a given root system $R$ in $V$ and a given multiplicity function $k$ he wrote down a trigonometric variant of the Dunkl operators (3.1) for $\xi \in V$:

$$(D_\xi^{(k)} f)(x) := (\partial_\xi f)(x) + \tfrac{1}{2} \sum_{\alpha \in R_+} k_\alpha \langle \alpha, \xi \rangle \cot(\tfrac{1}{2}\langle \alpha, x \rangle)\,(f(x) - f(s_\alpha x))$$
$$(x \in V,\ f \in C^\infty(V)). \qquad (4.7)$$

Now the operators $D_\xi^{(k)}$ will no longer commute, in general. However, Heckman showed that the operators $\sum_{\eta \in W\xi}(D_\eta^{(k)})^j$ ($\xi \in V$, $j = 0, 1, 2, \ldots$), when restricted to the $W$-invariant $C^\infty$ functions on $V$, coincide with differential operators which commute with each other and form a commutative algebra. This is the algebra looked for in Problem 2 of the previous subsection. The Jacobi polynomials $P_\lambda^{(k)}$ are the joint eigenfunctions of the operators in this algebra. By this approach, Heckman also obtained a quick existence proof for Opdam's shift operators.

Next Cherednik [3] made a slight but significant variation in Heckman's Dunkl type operators (4.7). He put

$$(\widetilde{D}_\xi^{(k)} f)(x) := (\partial_\xi f)(x) + \sum_{\alpha \in R_+} k_\alpha \langle \alpha, \xi \rangle \frac{1}{1 - e^{-\alpha}(x)}\,(f(x) - f(s_\alpha x))$$
$$- \tfrac{1}{2} \sum_{\alpha \in R_+} k_\alpha \langle \alpha, x \rangle\,f(x).$$

(Here I took the part of the right hand side on the second line from Opdam [13, p.86]; Cherednik is not very specific about this part of his formula.) Cherednik's operators have the nice property that they mutually commute, without the need of first restricting to $W$-invariant functions. On the other hand, they do not share the property $w\,D_\xi^{(k)}\,w^{-1} = D_{w\xi}^{(k)}$ of Heckman's operators. Anyhow, by means of Cherednik's operators one can draw the same conclusions as by Heckman's operators, and in a similar way. Moreover, a structure of graded Hecke algebra can be associated with Cherednik's operators.

# 5. Macdonald polynomials associated with $R$

## 5.1. Definition of Macdonald polynomials.
Let $0 < q < 1$. We keep the assumptions of §4.1 except that we replace the weight function $\delta_k$ in (4.1) by

$$\delta_{k,q}(x) := \prod_{\alpha \in R_+} \left| \frac{(e^{i\langle \alpha, x \rangle}; q)_\infty}{(q^{k_\alpha} e^{i\langle \alpha, x \rangle}; q)_\infty} \right|^2.$$

Then the *Macdonald polynomials* $P_\lambda^{(k,q)}$ were defined by Macdonald [23], [24] just as the Jacobi polynomials $P_\lambda^{(k)}$, but with the inner product in (4.2) replaced by

$$\langle f, g \rangle_{k,q} := \int_T f(x) \overline{g(x)} \delta_{k,q}(x) \, dx \quad (f, g \in \mathcal{A}).$$

In the case of root system $A_1$ the Macdonald polynomials coincide with the $q$-ultraspherical polynomials $x \mapsto C_n^{k,q}(\cos x)$. For any root system $R$, in the limit for $q \uparrow 1$, the Macdonald polynomial $P_\lambda^{(k,q)}$ tends to the corresponding Jacobi polynomial $P_\lambda^{(k)}$.

Macdonald gives some explicit $q$-difference operators of which the $P_\lambda^{(k,q)}$ are eigenfunctions. Although these operators, except for root system $A_n$ (where they were independently found by Ruijsenaars [30]) do not yet give a full commutative algebra of operators having the $P_\lambda^{(k,q)}$ as joint eigenfunctions, the additional parameter $q$ gives enough freedom such that already the eigenvalue of one such operator separates the elements of $P_+$ for generic $q$, by which a positive answer to question 1 in §4.2 can be given for the case of Macdonald polynomials. Taking limits for $q \uparrow 1$ then yields the same positive answer for the case of Jacobi polynomials. This is an alternative to Heckman's approach via Problem 2. Macdonald also gives conjectured explicit expressions for the squared norms $\langle P_\lambda^{(k,q)}, P_\lambda^{(k,q)} \rangle_{k,q}$.

## 5.2. Askey-Wilson polynomials for root system $BC_n$.
The author [18] introduced for the non-reduced root system $BC_n$ a class of polynomials having two more parameters than Macdonald's class for $BC_n$. This extended class reduces for $n = 1$ to the Askey-Wilson polynomials [2]. In [18] only one explicit $q$-difference operator was given having the $BC_n$-polynomials as eigenfunctions, but this was sufficient for establishing orthogonality. Later, van Diejen [7] gave explicit expressions for the generators af a full commutative algebra of operators having the $BC_n$ polynomials as joint eigenfunctions.

## 5.3. Cherednik's approach to Macdonald polynomials.
Cherednik [4], [5] succeeded to give positive answers to questions 2 and 3 in §4.2. In the context of certain representations of affine Hecke algebras he could realize a commutative algebra of operators which have the Macdonald polynomials as joint eigenfunctions. In the same context he could realize $q$-analogues of Opdam's shift operators and next, by the same technique as in Opdam, prove Macdonald's conjectures in the $q$-case.

It is beyond the scope of this short survey to explain Cherednik's approach in any detail. In May 1994 I. G. Macdonald delivered some very helpful lectures in Leiden in order to explain Cherednik's approach. Let me here only give a few indications. Just as a Hecke algebra is a deformation of the group algebra of a Weyl group, an affine Hecke algebra (cf. [20]) deforms the group algebra of an affine Weyl group. If $R$ is an irreducible root system in $V$ with Weyl group $W$ then the (extended) *affine Weyl group* is the semidirect product $\widetilde{W} := W \ltimes P\check{\,}$, where the *dual weight lattice* $P\check{\,}$ is defined as $P\check{\,} := \{\lambda \in V \mid \langle \lambda, \alpha \rangle \in \mathbb{Z}\}$, an abelian group under addition. Then $\widetilde{W}$ acts as a group of motions on $V$, with $P\check{\,}$ acting as a group of translations. The group $\widetilde{W}$ also acts on $\mathcal{A}$, with $W$ acting as before and with the action of $P\check{\,}$ still depending on a parameter $q$.

The *affine Hecke algebra* $H$ can be defined in terms of generators and relations which still depend on the values of a $W$-invariant function $\alpha \mapsto t_\alpha$ on $R$. Corresponding to a choice of $R_+$ we can define $P_+\check{\,}$. The embedding of $P_+\check{\,}$ in $H$ then generates a commutative subalgebra $\mathcal{Y}$ of $H$.

For given $q$ we can use the action of $\widetilde{W}$ on $\mathcal{A}$ in order to define an action of $H$ on $\mathcal{A}$, by specifying the action for a set of generators of $H$ (Demazure operators). This action depends on $q$ and the $t_\alpha$. Put $t_\alpha = q^{-k_\alpha/2}$, where $k \colon R \to [0, \infty)$ is a multiplicity function. Then Cherednik proves that the Macdonald polynomials $P_\lambda^{(k,q)}$ are the joint eigenfunctions of the $W$-invariant elements in the commutative algebra $\mathcal{Y}$. This answers question 2 in §4.2.

At the moment it is still an open problem to extend Cherednik's approach to the $BC_n$ polynomials. As van Diejen [7] already answered question 2 in §4.2 in a positive way for this case, it would be nice to complement van Diejen's constructive approach with the deep conceptual approach via affine Hecke algebras.